\documentclass{article}
\usepackage{amssymb,amsmath,amsthm}
\usepackage{graphicx,graphics,epsfig}
\newtheorem{theorem}{Theorem}

\newtheorem{lemma}{Lemma}
\newtheorem{proposition}{Proposition}

\newtheorem{remark}{Remark}
\date{}
\numberwithin{equation}{section} \numberwithin{theorem}{section}
\numberwithin{lemma}{section} \numberwithin{corollary}{section}
\numberwithin{remark}{section} \numberwithin{proposition}{section}
\numberwithin{definition}{section}
\newcommand{\n}{\noindent}
\newcommand{\vs}{\vskip}
\begin{document}

\title{A class of free boundary problems with Neuman boundary condition}

\vs 0.5cm
\author{ABDESLEM LYAGHFOURI\\
American University of Ras Al Khaimah,\\ Department of Mathematics and Natural
Sciences\\
Ras Al Khaimah, UAE\\
ABDERACHID SAADI\\
Mohamed Boudiaf University,\\ Department of Mathematics, Msila, Algeria\\
and\\
Laboratoire d'Equations aux D\'{e}riv\'{e}es Partielles Non Lin\'{e}aires\\
et Histoire des Math\'{e}matiques\\
ENS, Kouba, Algeria} \maketitle

\begin{abstract}
In this work, we investigate the continuity of the free boundary in a class of elliptic problems, with
Neuman boundary condition. The main idea is a change of variable
that allows us to reduce the problem to the one studied in  \cite{[S]}.
\end{abstract}
\vs 0.3cm
\n Mathematics subject classification: 35J15, 35R35
\vs 0.2cm
\n Keywords: Neuman boundary condition, Free boundary, Change of variable

\section{Statement of the problem  and preliminary results}\label{1}
\vs 0.5cm  In \cite{[S]}, the second author studied the following class of problems:
\begin{equation*}(P_0)
\begin{cases}
& \text{Find } (u, \chi) \in  H^{1}(\Omega)\times L^\infty
(\Omega)
\text{ such that}:\\
& (i)\quad  u\geq 0, \quad 0\leq \chi\leq 1 ,
\quad u(\chi -1 ) = 0 \,\,\text{ a.e.  in } \Omega\\
& (ii)\quad  u= 0  \quad \text{ on } \Gamma_2 \\
& (iii)\quad \displaystyle{\int_\Omega }\big( a(x) \nabla
u + \chi h(x) \big) .\nabla\xi dX \,\leq\, \int_{\Gamma_3}\beta(x,\varphi-u)\xi d\sigma(x)\\
&\hskip 1.7cm \forall \xi \in H^{1}(\Omega),\quad\xi \geq 0 \text{
on } \Gamma_2
\end{cases}
\end{equation*}
\n
where $\Omega=\{(x_1,x_2)\in \mathbb{R}^2\,\,/\,\,x_1\in(a_0,b_0),\,\, d_0<x_2<\gamma(x_1)\}$,
with $\gamma\in C^{0,1}(a_0,b_0)$,
$\Gamma_2\cup \Gamma_3=\{(x_1,\gamma(x_1))\,\,/\,\,x_1\in(a_0,b_0)\}$,
$\Gamma_2\cap \Gamma_3=\emptyset$, $\Gamma_3$
is a nonempty connected and relatively open subset of $\partial\Omega$,
$a=[a_{ij}]$ is a $2\times2$ matrix and $h$ is a nonnegative function.

\n In \cite{[C]}, \cite{[ChaL2]}, and \cite{[S]}, the monotonicity  of
$\chi$ with respect to the variable $x_2$, has allowed the authors to define
the free boundary $\partial \{u>0\} \cap\Omega$ as a graph of a function 
$\Phi(x_1)$.
Moreover, under suitable assumptions, it was proven that $\Phi$ is continuous
both for Dirichlet and Neuman conditions.

\vs0.2cm\n In this paper, we consider a more general class
of free boundary  problems  in the same spirit of \cite{[ChaL3]},
namely we replace the real valued function $h$ by a vector function
$H$:

\begin{equation*}(P)
\begin{cases}
& \text{Find } (u, \chi) \in  H^{1}(\Omega)\times L^\infty (\Omega)
\text{ such that}:\\
& (i)\quad  u\geq 0, \quad 0\leq \chi\leq 1 ,
\quad u(1-\chi) = 0 \,\,\text{ a.e.  in } \Omega\\
& (ii)\quad  u= 0  \quad \text{ on } \Gamma_2 \\
& (iii)\quad \displaystyle{\int_\Omega }\big( a(x) \nabla
u + \chi H(x)\big) .\nabla\xi dx \,\leq\, \int_{\Gamma_3}\beta(x,\varphi-u)\xi d\sigma(x)\\
&\hskip 1.7cm \forall \xi \in H^{1}(\Omega),\quad\xi \geq 0 \text{
on } \Gamma_2
\end{cases}
\end{equation*}

\n where $\Omega$ is a bounded domain of $\mathbb{R}^2$ whose
boundary $\partial\Omega$ is of class $C^1$,
$\Gamma_2$ and $\Gamma_3$ are disjoint nonempty subsets of $\partial\Omega$,
with $\Gamma_3$ connected and relatively open in $\partial\Omega$.

\n $a=[a_{ij}]$ is a $2\times2$ matrix that satisfies for two positive constants
$\lambda$ and $\Lambda$
\begin{eqnarray}
 & \vert a_{ij}(x)\vert \leq \Lambda,
 \quad\text{  for a.e. }x\in \Omega, \quad \forall i, j=1, 2\\
&{ a}(x)\xi.\xi \geq \lambda\vert \xi\vert^2\quad\forall\xi\in
\mathbb{R}^{2},\quad\text{  for a.e. }x\in \Omega,
\end{eqnarray}

\n $H=(H_1,H_2)$ is a $C^1(\overline{\Omega})$ vector function, that satisfies for
some positive constants $\overline{H}>\underline{H}$:
\begin{eqnarray}
&  |H_1(x)|\leq \bar H \quad\text{ in } \Omega  \\
& \underline{H} \leq H_2(x)\leq \overline{H} \quad\text{ in } \Omega  \\
& \text{div}(H(x))\geq 0\quad\text{ in }\Omega\\
& H(x).\nu>0\quad \text{ on } \Gamma_3
\end{eqnarray}

\begin{eqnarray}
&& \beta(x,.)\quad\text{ is continuous for a.e. }x\in \Gamma_3 \\
&& \beta(x,0)=0\quad\text{for a.e. }x\in \Gamma_3 \\
&& \beta(x,.)\quad\text{ is non-decreasing for a.e. }x\in \Gamma_3
\end{eqnarray}

\n Many free boundary problems belongs to the above class of problems.
For example the dam problem with Neuman boundary condition on the
reservoirs bottoms (see \cite{[ChiL1]}, \cite{[ChiL2]}, \cite{[L1]},
\cite{[L2]}, \cite{[L3]}, \cite{[L4]}). Another problem arises from
the thermoelectrical modelling of aluminum electrolytic cells (see \cite{[BMQ]}).

\vs0.2cm\n In these problems we investigate the free boundary
$\partial[u>0]\cap\Omega$ that separates two different regions.
In the dam problem, it separates wet and non wet parts of
the porous medium. In the aluminium electrolysis problem, it
separates liquid and solid aluminium.

\begin{remark}\label{r1.1} Under assumptions (1.1)-(1.4) and (1.7)-(1.9), we can prove
the existence of a solution for the problem $(P)$ as in \cite{[ChiL1]}.
For a more general situation, we refer to \cite{[L1]}.
\end{remark}

\vs0.3cm \n We begin by the following proposition that can be obtained as in \cite{[ChaL3]}.
\begin{proposition}\label{p1.1}
For any solution $(u,\chi)$ of $(P)$, we have:
\begin{eqnarray*}
&& i)\quad \text{div}(a(x)\nabla u)=-\text{div}(\chi H(x)) \quad\text{in}\quad \mathcal{D}'(\Omega).\\
&&ii)\quad \text{div}(\chi H(x))-\chi_{\{u>0\}}\text{div}(H(x))\leq0\quad\text{in}\quad\mathcal{D}'(\Omega).
\end{eqnarray*}
\end{proposition}

\begin{remark}\label{r1.2} As a consequence of Proposition 1.1 i),
we obtain (see \cite{[GT]}):

\vs 0.2cm\n $i)$ $u\in C_{loc}^{0,\alpha}( \Omega)$ for some
$\alpha\in(0,1)$. In particular $u$ is continuous in $\Omega\cup\Gamma_2$
and the set $\{u>0\}$ is open.

\vs 0.2cm \n $ii)$ If $a\in C_{loc}^{0,\alpha}( \Omega)$
$(0<\alpha<1)$, then we have $u\in C_{loc}^{1,\alpha} (\{u>0\})$.
\end{remark}

\vs 0,5cm\n Following \cite{[ChaL1]}, we introduce for each $h\in \pi_{x_{2}}(\Omega)$ and $w\in
\pi_{x_1}(\Omega\cap\{x_{2}=h\})$, the following differential equation:

$$(E(w,h))\left\{\begin{array}{l}
  X' (t,w,h)=  H(X(t,w,h))\\
  X(0,w,h)=(w,h)\\
  \end{array}\right.$$

\n It is well known that this equation has a maximal solution $X(.,w,h)$ defined on a
maximal interval $(\alpha_{-}(w,h), \alpha_{+}(w,h))$ and continuous on the open set:
$$\{(t,w,h):~\alpha_{-}(w,h)< t < \alpha_{+}(w,h),
 h\in \pi_{x_2}(\Omega) , w\in \pi_{x_1}(\Omega\cap\{x_2=h\})\}$$

\n Moreover due to (1.4), we have:
\[X(\alpha_{-}(w,h),w,h) \in \partial\Omega\cap\{x_2<h\} \quad
\text{and}\quad X(\alpha_{+}(w,h),w,h) \in \partial\Omega\cap\{x_2>h\}\]

\n In the sequel, we will denote the functions $X(t,w,h), \alpha_{-}(w,h)$, and
$\alpha_{+}(w,h)$ respectively by $X(t,w),\alpha_{-}(w),$ and $\alpha_{+}(w)$.

\vs 0.3cm \n The function $\alpha_{-}$ (resp. $\alpha_{+}$) is upper
(resp. lower) semi-continuous. The next result gives more regularity
for $\alpha_{+}$.

\begin{theorem}\label{1.1} For every $h\in\pi_{x_2}(\Omega)$, $\alpha_{+}$ is
continuously differentiable at each $w_0\in\pi_{x_1}(\Omega\cap \{x_2=h\})$
such that $x_0=X(\alpha_{+} (w_0),w_0)\in \Gamma_3$.
\end{theorem}

\n \emph{Proof.} Let $h$ and $w_0$ as in the theorem.
Since $\partial\Omega$ is a $C^1$ curve, there exists
an open set $U\subset \mathbb{R}^2$  that contains $x_0$
and a $C^1-$diffeomorphism $\Upsilon=(\Upsilon_1\Upsilon_2):~U~\rightarrow~B_1$ such that
\begin{equation}\label{1.10}
\Upsilon(U\cap\Omega)=B_1\cap\{y_2>0\}\quad\text{and}\quad \Upsilon(U\cap\partial\Omega)=B_1\cap\{y_2=0\},
\end{equation}
where $B_1$ is the unit ball.

\vs 0.2 cm\n Let $x_0^-\in (U\cap\partial\Omega)\setminus\{x_0\}$ such that $(x_0^--x_0).\tau(x_0)<0$,
where $\tau(x_0)$ is the unit tangent vector to $\partial\Omega$ at $x_0$.

\vs 0.2 cm\n Since $H\in C^1(\overline{\Omega})$, there exists an open set
$\Omega^*$ and an extension $H^*$ of $H$ such that $\bar\Omega\subset \Omega^*$ and $H^*\in
C^1( \Omega^*)$. Then we consider the unique maximal solution $Z(t)$ of the differential equation:
 \[\left\{\begin{array}{l}
  Z' (t)=  H^*(Z(t))\\
  Z(0)=x_0^-\\
  \end{array}\right.\]
which is defined on a maximal open interval $(\gamma,\delta)$.

\vs 0.2 cm\n Taking into account (1.6), we can see that $Z(t)\in \Omega$ for all $t\in(\gamma,0)$.
Now if we assume that $h$ is close enough to $x_{02}$, and denote by $t_h$ the real number for 
which the curve $Z(t)$ intersects the line $x_2=h$, then there exists $w_0^-\in \pi_{x_1}(\Omega\cap
\{x_2=h\})$  such that $Z(t_h)=(w_0^-,h)$. Moreover, it is easy to see that
\[\left\{\begin{array}{l}
  X(t)=Z(t_h-t)\\
  X(0)=(w_0^-,h)\\
  \end{array}\right.\]

\n Since $(x_0^--x_0).\tau(x_0)<0$, we necessarily have  $w_0^-<w_0$.
Furthermore, for each $w_0^-<w<w_0$, the curve $X(t,w)$
is located between the curves $X(t,w_0)$ and $X(t,w_0^-)$.
Therefore we have
\begin{equation}\label{1.11}
X( \alpha_{+}(w),w)\in U\cap\partial\Omega\quad \forall w\in
(w_0^-,w_0)
\end{equation}

\vs 0.2 cm\n Let now $x_0^+\in (U\cap\partial\Omega)\setminus\{x_0\}$ such that $(x_0^+-x_0).\tau(x_0)>0$.
Arguing as above, we can prove that there exists $w_0^+\in \pi_{x_1}(\Omega\cap
\{x_2=h\})$  such that
\begin{equation}\label{1.12}
X( \alpha_{+}(w),w)\in U\cap\partial\Omega\quad \forall w\in
(w_0,w_0^+)
\end{equation}

\n Taking into account (1.10)-(1.12), we see that there exists $\eta>0$ small enough such that
\begin{equation}\label{1.13}
\Upsilon_2(X( \alpha_{+}(w),w))=0\quad \forall w\in
(w_0-\eta,w_0+\eta)
\end{equation}

\vs 0.2cm\n For each $\omega\in \pi_{x_1}(\Omega^*\cap\{x_2=h\})$, let $X^*(t,w)$
be the unique maximal solution of the differential equation:
 \[(E(w,h))\left\{\begin{array}{l}
  X' (t,w)=  H^*(X(t,w))\\
  X(0,w)=(w,h)\\
  \end{array}\right.\]
$X^*(t,w)$ is defined on the interval $[\alpha_{-}^* (w), \alpha_{+}^*(w)]$, and we obviously have
$X^*_{\vert_{(\alpha_{-}(w), \alpha_{+}(w))}}=X$ .
Moreover, we have $\alpha_{-}^* (w)<\alpha_{-}(w)$ and $\alpha_{+}(w)<\alpha_{+}^*(w)$.

\vs 0.2cm\n Let $D^*=\{ (t,w) \,/\, w\in (w_0-\eta,w_0+\eta),\, t\in (\alpha_{-}^* (w), \alpha_{+}^*(w))\}$.
Since $X^*\in C^1(D^*)$ and $\Upsilon_2\in C^1(U)$, the function $F^*=\Upsilon_2 oX^*$ is in $C^1(D^*)$.
In addition, $F^*$ is an extension of $F=\Upsilon_2 oX$ to $D^*$ and we have
\begin{eqnarray*}{{\partial F^*}\over{\partial t}}(t,w)
&=& \nabla \Upsilon_2(X^*(t,w)).X'^*(t,w)\\
&=& \nabla \Upsilon_2(X^*(t,w)).H^*(X^*(t,w))
\end{eqnarray*}

\n In particular, we obtain from (1.6) and (1.13)
\[{{\partial F^*}\over{\partial t}}(\alpha_{+} (w_0), w_0)
=\nabla \Upsilon_2(X(\alpha_{+} (w_0), w_0)) . H(X(\alpha_{+} (w_0), w_0))\neq0
\]
Therefore by the implicit function theorem, there
exists $\delta \in(0,\eta)$ and a unique function $f : (w_0
-\delta , w_0+\delta) \rightarrow {I\!\! R}$ such that
\begin{eqnarray*}
 F^* (t,\omega) &=& 0 \quad\text{iff}\quad t=f(\omega)\\
f(w_0) &=&\alpha_{+}(w_0) \quad\text{ and }\quad f\in
C^1 (w_0-\delta,w_0+\delta).\end{eqnarray*} Since $F^*
(\alpha_{+} (w),w)= F (\alpha_{+} (w),w)=0$,
we have $\alpha_{+}(w)= f(w)$ for all $w\in(w_0-\delta,w_0+\delta)$.
We conclude that $\alpha_{+}\in C^1(\pi_{x_1}(\Omega\cap\{x_2=h\}))$.
\qed

\vs 0.3cm\n Following \cite{[ChaL1]}, we define
for each $h\in\pi_{x_2}(\Omega)$, the set:
$$D_{h}=\{(t,w): w\in\pi_{x_1}(\Omega\cap\{x_2=
h\}), t\in (\alpha_{-}(w), \alpha_{+}(w))\}$$

\n and the $C^1$ mapping:
$$\begin{array}{llll}
T_{h} :& D_{h} & \longrightarrow & T_{h}(D_{h})\\
     & (t,w) &\longmapsto& T_{h}(t,w)= X(t,w)\\
  \end{array}$$

\n  whose Jacobian determinant is denoted by $Y_{h}(t,w)$.

\vs 0,3cm\n The next proposition can be established as in \cite{[ChaL1]}:
\begin{proposition}\label{p1.2}

\n $i)$ $T_{h}$ is a $C^1-$diffeomorphism.

\n $ii)$ $\displaystyle \frac{\partial Y_{h}}{\partial t}(t,w)=
Y_{h}(t,w).div(H(X(t,w))).$

\n $iii)$ $Y_{h}(t,w)=-H_2(w,h)exp\left[\displaystyle\int_{0}^{t}
div(H(X(s,w))ds\right].$
\end{proposition}

\vs0.2cm\n In Section 3, we will use the $C^1-$diffeomorphism $T_h$ as a
change of variable to transform the problem $(P)$ to a problem of type $(P_0)$.
As a consequence, we obtain from \cite{[S]} that the free boundary is locally
represented by graphs of a family of continuous functions.

\section{Parametrization of the free boundary}\label{2}

\vs 0.5cm\n For each $h\in\pi_{x_2}(\Omega)$ and any function $f$ defined
in $\Omega$, we shall denote by $\widetilde{f}$ the function $f\circ T_{h}$.

\n The first result of this section is a monotonicity property of $\widetilde{\chi}$
with respect to $t$, which translates into the fact that $\chi$ is non-increasing 
along the orbits of the differential equation $E(w,h)$. For the proof we refer to the
one of Theorem 2.1 in \cite{[ChaL3]}.
\begin{proposition}\label{p2.1}
Let $(u,\chi)$ be a solution of $(P)$. Then we have for each $h\in
\pi_{x_2}(\Omega)$:
\[\frac{\partial\widetilde{\chi}}{\partial t} \leq0 \;\;\;in\;\;\;\mathcal{D}'(D_{h})\]
\end{proposition}

\vs 0,5cm\n The next proposition is a consequence of the monotonicity of $\widetilde{\chi}$ and
the continuity of $\widetilde{u}$. For the proof we refer to the
one of Proposition 3.1 in \cite{[ChaL3]}

\begin{proposition} \label{p2.2} Let $(u,\chi)$ be a solution of $(P)$ and
$(t_0,w_0)\in D_h$.

\vs 0,2cm \n
 $i)$\quad If $\widetilde{u}(t_0,w_0)>0$, then there exists $\epsilon>0$ such that:
    $$\widetilde{u}(t,w)>0\quad \forall (t,w)\in
    \mathcal{C}_{\epsilon}=\{(t,w)\in D_h:
    |w-w_{0}|<\epsilon,~ t<t_0+\epsilon\}$$

\vs 0,2cm \n $ii)$\quad If $\widetilde{u}(t_0,w_0)=0$, then:
    $$\widetilde{u}(t,w_0)=0, \quad\forall t\geq t_0$$
\end{proposition}

\vs 0,5cm \n Thanks to Proposition 2.2, we can define for each $h\in\pi_{x_2}(\Omega)$,
the following function in $\pi_{x_1}(\Omega\cap\{x_2=h\})$:
\begin{equation*}
    \Phi_{h}(w)=\left\{\begin{array}{lll}
 \sup \{t:(t,w)\in D_{h}:\widetilde{u}(t,w)>0\} & : &
 \text{if  this  set  is  not  empty} \\
  \alpha_{-}(w) & : & \text{otherwise} \\
\end{array}\right.
\end{equation*}

\vs 0,2cm\n Arguing as in \cite{[ChaL1]}, we can see that $\Phi_h$
is well defined and satisfies

\vs 0,2cm
\begin{proposition}\label{p2.3} $\Phi_h$ is lower
semi-continuous on $\pi_{x_1} (\Omega\cap\{x_2=h\})$ and
$$\{\widetilde{u}>0\}\cap D_h=\{t<\Phi_h(w)\}$$
\end{proposition}

\begin{remark}\label{r1.3} If the functions $\Phi_h$ are smooth,
then the family of functions $\{\Phi_h\}$ is a local parametrization of
the free boundary $\partial \{u>0\} \cap\Omega$.
\end{remark}

\vs 0,5cm\n The next result gives a description of $\chi$ in the interior
of the set $\{u=0\}$.

\begin{theorem} \label{th2.1} Let $(u,\chi)$ be a solution of $(P)$,
$(x_{01},x_{02})=T_{h}(t_0,w_0) \in T_{h}(D_h)$,
$B_r(t_0,w_0)$ the ball of center $(t_0,w_0)$ and radius $r$,
$Z_0=\big((t_0,\infty)\times (w_0-r,w_0+r)\big)\cap D_h$ and
$C_r=Z_0 \cup B_{r}(t_0,w_0)$.

\n If $\widetilde{u}=0$ in $B_r(t_0,w_0)\subset D_h$, then we have
$\widetilde{u}=0$ in $C_r$. Moreover
\begin{enumerate}
\item If $\overline{T_{h} (Z_{0})}
  \cap \Gamma_{3} =\emptyset,$ then $\widetilde{\chi}=0$ in $C_r$.
\item If $\overline{T_{h} (Z_{0})} \cap \Gamma_{2} = \emptyset,$ then:
\[
\widetilde{\chi}(t,w)=\displaystyle\frac{Y_{h}(\alpha
_{+}(w),w)} {Y_{h}(t,w)}\frac {\beta (., \varphi(.))}
{H.\nu}(X(\alpha_{+}(w),w)).
\]
\end{enumerate}
\end{theorem}

\n To prove the theorem, we need two lemmas.

\vs 0.3cm
\begin{lemma}\label{lem2.1}
For each $x_0\in \Gamma_3$, there exists $\eta>0$ small enough and a $C^1$ function
$\sigma$ such that one of the following
situations holds
\begin{eqnarray*}
i) &\quad   \Gamma_3\cap B(x_0,\eta)\subset\{(x_1,\sigma(x_1))~\} \\
ii)&\quad  \Gamma_3\cap B(x_0,\eta)\subset\{(\sigma(x_2),x_2)~\} \\
\end{eqnarray*}
\end{lemma}

\n\emph{Proof.} Since $\Gamma_3$ is a $C^1-$curve, there exists
an open set $U\subset \mathbb{R}^2$  that contains $x_0=(x_{01},x_{02})$
and a $C^1-$diffeomorphism $\Upsilon:~U~\rightarrow~B_1$ such that
$\Upsilon(U\cap\Omega)=B_1\cap\{y_2>0\}$ and $\Upsilon(U\cap\Gamma_3)=B_1\cap\{y_2=0\}$.

\vs0.2cm\n If $\Upsilon=(\Upsilon_1\Upsilon_2)$, then we have:
\[\Upsilon_2(x)=0\quad \forall x\in U\cap\Gamma_3\]

\vs0.2cm\n Due to (1.6), we have $\nabla\Upsilon_2(x_0)\neq0$.
Therefore either $\displaystyle{{{\partial \Upsilon_2}\over{\partial x_1}}(x_0)\neq0}$,
or $\displaystyle{{{\partial \Upsilon_2}\over{\partial x_2}}(x_0)\neq0}$.

\vs0.2cm\n Assume for example that we have $\displaystyle{{{\partial \Upsilon_2}\over{\partial x_2}}(x_0)\neq0}$.
Then by the implicit function theorem, there
exists $\delta>0$ small enough and a unique $C^1-$function $\sigma : (x_{01}
-\delta , x_{01}+\delta) \rightarrow {I\!\! R}$ such that
\begin{eqnarray*}
&&\Upsilon_2(x_1,x_2)= 0 \quad\text{iff} \quad x_2=\sigma(x_1)\\
&&\quad\text{for all } x_1\in(x_{01}-\delta,x_{01}+\delta).
\end{eqnarray*}

\n So $i)$ holds.

\vs 0.2cm\n If $\displaystyle{{{\partial \Upsilon_2}\over{\partial x_1}}(x_0)\neq0}$,
then we can show in a same fashion that $ii)$ holds.

\qed

\begin{lemma}\label{lem2.2}
Let $w_1, w_2\in \pi_{x_1}(\Omega\cap\{x_2=h\})$ such that $w_1<w_2$
and $T_h(\alpha_+(w_i),w_i)\in \Gamma_3$, for $i=1, 2$. Then we have:
\begin{eqnarray*}
&&\int_{Z}\big( \texttt{a}(x) \nabla\tilde{u} + \tilde{\chi} \texttt{h}(t,w)e_t\big) .\nabla\xi dtdw \,=\,
\int_{\tilde{\Gamma}_3}\lambda(.,\tilde{\varphi}-\tilde{u})\xi d\tilde{\sigma}\\
&&\qquad\forall \xi \in H^1(Z),\quad\xi=0 ~~\text{ on }~~ \partial Z\cap D_h
\end{eqnarray*}
where
\begin{eqnarray*}
&& Z=\{(t,w):~w_1<w<w_2 ~\text{ and }~ h< t < \alpha_{+}(w)\}\\
&& \tilde{\Gamma}_3=\{(\alpha_{+}(w),w):~w_1<w<w_2\}
\end{eqnarray*}

\begin{eqnarray*}
&& \lambda((t,w),z)=\mu(w)\beta(T_h(t,w),z)\\
&&\mu(w)={{|Y_h|(\alpha_{+}(w),w)}\over
\sqrt{1+\alpha_{+}^{\prime 2}(w)}(H.\nu)(T_h(\alpha_{+}(w),w)) }\\
&& \texttt{h}(t,w) = |Y_h(t,w)|, \qquad e_t=(1,0) \\
&&  \texttt{a}(t,w)= |Y_h(t,w)|  ^{t}P(t,w).a(X(t,w)).P(t,w) \\
&&\hbox{with } \quad P= (^{t}\mathcal{J}T_h)^{-1} = \displaystyle{1\over Y_h(t,w)}
    \left(%
\begin{array}{cc}
 \displaystyle{\partial X_2\over \partial\omega} (t,w)
 & -H_2(X(t,w))\\
-\displaystyle{\partial X_1\over \partial\omega} (t,w) &
 H_1(X(t,w))\\
\end{array}%
\right).
 \end{eqnarray*}
\end{lemma}

\vs0,3cm\n\emph{Proof.} Let $\xi\in H^1(Z)$ such that $\xi=0$ on $\partial Z\cap D_h$.
Then  $\pm\xi o T_h^{-1} \chi( T_h(Z))$ are test functions for $(P)$ and we have

\begin{equation}\label{2.1}
\int_{T_h(Z)}( a(x)\nabla u + \chi H(x)). \nabla(\xi o
 T_h^{-1}) dx=\int_{\Gamma_3\cap T_h(\partial Z)}\beta(x,\varphi-u)\xi oT_h^{-1} d\sigma(x)
 \end{equation}

\n The left hand side of (2.1) can be written using the change of variable $T_h$
(see \cite{[ChaL3]}) as

\begin{equation}\label{2.2}
\int_{D_h}( \texttt{a}(t,\omega)\nabla (u oT_h) + \chi oT_h . \texttt{h}(t,\omega) e_t). \nabla \xi
  dt d\omega
\end{equation}
where the matrix $\texttt{a}$ and the function $\texttt{h}$ are
given by

\begin{eqnarray*}
&& \texttt{h}(t,\omega) = |Y_h(t,\omega)|, \qquad e_t=(1,0) \\
&&  \texttt{a}(t,\omega)= |Y_h(t,\omega)|  ^{t}P(t,\omega).a(X(t,\omega)).P(t,\omega) \\
&&\hbox{with } \quad P= (^{t}\mathcal{J}T_h)^{-1} = \displaystyle{1\over Y_h(t,\omega)}
    \left(%
\begin{array}{cc}
 \displaystyle{\partial X_2\over \partial\omega} (t,\omega)
 & -H_2(X(t,\omega))

 \\\\
-\displaystyle{\partial X_1\over \partial\omega} (t,\omega) &
 H_1(X(t,\omega))\\
\end{array}%
\right).
 \end{eqnarray*}

\n To handle the right hand side of (2.1), we first observe that
\begin{equation}\label{2.3}
\{T_h(\alpha^+(w),w), w_1<w<w_2\}=\Gamma_3\cap  T_h(\partial Z)
\end{equation}

\n Shrinking if necessary, we can assume by Lemma 2.1, that there
exists a $C^1-$function $\sigma$ such that one of the following
situations holds
\begin{eqnarray*}
i)&\quad   \sigma(X_1 (\alpha_{+}(w),w))= X_2(\alpha_{+}(w),w)\quad \forall w\in
(w_1,w_2)  , \\
ii)&\quad   \sigma(X_2(\alpha_{+}(w),w))= X_1 (
\alpha_{+}(w),w)\quad \forall \omega\in (w_1,w_2).
\end{eqnarray*}
Assume for example that $i)$ holds. The case $ii)$ can
be treated in the same way. Since $x_1\rightarrow (x_1,\sigma(x_1))$ is a
$C^1-$parametrization of $\Gamma_3\cap\partial (T_h (Z))$, the
integral in the right hand side of (2.4) can be written as

\begin{eqnarray}\label{2.4}
&&\int_{\Gamma_3\cap T_h(\partial Z)}\beta(x,\varphi-u)\xi oT_h^{-1} d\sigma(x)\nonumber\\
&&=\int_{\pi_{x_1}(\Gamma_3\cap\partial
(T_h(Z))}\beta((x_1,\sigma (x_1)), \varphi
(x,\sigma(x)))\xi oT_h^{-1}(x_1,\sigma (x_1))\sqrt{1\!+\!\sigma^{\prime 2 }(x_1) }dx_1\nonumber\\
\end{eqnarray}

\n Now observe that $(x_1,\sigma(x_1))= T_h (\alpha_{+}(w),w)$
for $w\in (w_1,w_2)$, and let $\theta(w) =
x_1= T_h^1( \alpha_{+}(w),w)$. Then $\theta$ is a $C^1-$function
and $\theta'(w) = \alpha_{+}'(w) H_1(X(\alpha_{+}(w),w)) +{\partial X_1\over
\partial w} $. Using Theorem 1.1 and arguing as in \cite{[ChaL1]}, we can show via
implicit differentiation that
$$\alpha_{+}'(\omega)={\sigma' (X_{1}(\alpha_{+}(w),w)){\partial X_1
/\partial w}(\alpha_{+}(w),w) - {\partial X_2/\partial w}(\alpha_{+}(w),w)\over H_2(
X(\alpha_{+}(\omega),w)) - \sigma'(X_{1}(\alpha_{+}(w),w))H_1(X(\alpha_{+}(w),w)) }$$
which leads to
\begin{eqnarray*}
&&\theta'(w) = { -Y_h(\alpha_{+}(w),w)\over H_2(
X(w_{+}(w),w))-\sigma'(X_{1}(\alpha_{+}(w),w))H_1(X(\alpha_{+}(w),w)) }\\
&& = { |Y_h|(\alpha_{+}(w),w)(1+\sigma^{\prime 2 }(x_1))^{-1/2}\over H(X(\alpha_{+}(w),w),e).\nu
(X(\alpha_{+}(w),w)) }\end{eqnarray*}
where $\nu(x)={( -\sigma'(x_1),1)\over \sqrt{1+\sigma^{\prime 2 }(x_1)}}$ is the outward unit
normal to $\Gamma_3$.

\n Lastly we apply the change of variable $\theta$ to (2.4) to show that

\begin{eqnarray}\label{2.5}
&&\int_{\Gamma_3\cap T_h(\partial Z)}\beta(x,\varphi-u)\xi oT_h^{-1} d\sigma(x)\nonumber\\
&&=\int_{w_1}^{w_2} {\beta((T_h(\alpha_{+}(w),w))), \varphi
(T_h(\alpha_{+}(w),w)))  |Y_h|(\alpha_{+}(w),w)\over
H(T_h(\alpha_{+}(w),w)).\nu(T_h(\alpha_{+}(w),w)) }\xi(\alpha_{+}(w),w)) dw\nonumber\\
&&=\int_{w_1}^{w_2} {{\beta((T_h(\alpha_{+}(w),w))), \varphi
(T_h(\alpha_{+}(w),w)))|Y_h|(\alpha_{+}(w),w)}\over
\sqrt{1+\alpha_{+}^{\prime 2}(w)}H(T_h(\alpha_{+}(w),w)).\nu(T_h(\alpha_{+}(w),w)) }
\xi(\alpha_{+}(w),w)) d\sigma(w)\nonumber\\
&&=\int_{\widetilde{\Gamma}_3} \lambda((\alpha^+(w),w),\widetilde{\varphi}-\widetilde{u})\xi d\sigma(w)
\end{eqnarray}

\n Combining (2.1), (2.2) and (2.5), the result follows.
\qed

\vs0,3cm\n\emph{Proof of Theorem 2.1.} We first observe that $\widetilde{u}=0$ in $C_r$
and that statements 1) can be established as in \cite{[ChaL3]}.

\vs0,2cm\n Next we assume that $\overline{T_{h} (Z_0)} \cap \Gamma_{2} = \emptyset$.

\vs0,2cm\n From Lemma 2.2 and Proposition 2.4 of \cite{[S]}, we obtain for all $(t,w)$ in $C_r$:
\begin{eqnarray*}
\widetilde{\chi}(t,w)&=&\frac {\lambda((\alpha_{+}(w),w), \widetilde{\varphi}(\alpha_{+}(w),w))}
{\texttt{h}(t,w)\nu_2(\alpha_{+}(w),w)}\\
&=&\frac {{{|Y_h|(\alpha_{+}(w),w)}\over
\sqrt{1+\alpha_{+}^{\prime 2}(w)}H(T_h(\alpha_{+}(w),w)).\nu(T_h(\alpha_{+}(w),w)) }\beta(X(\alpha_{+}(w),w),\varphi(X(\alpha_{+}(w),w)))}
{|Y_h(t,w)|.\nu_2(\alpha_{+}(w),w)}\\
&=&{{|Y_h|(\alpha_{+}(w),w)}\over{|Y_h(t,w)|}}.{{\beta(.,\varphi)}\over{H.\nu}}(X(\alpha_{+}(w),w))
\end{eqnarray*}

\n Thus the result follows.
\qed

\vs 0,5cm

\section{Continuity of the Free Boundary}\label{3}

\vs 0,5cm\n In this section, we assume that:
\begin{eqnarray}\label{e3.1-3.4}
&& H\in C^{1,1}_{loc}( \Omega)\\
&&a\in C_{loc}^{0,\alpha}(\Omega\cup\Gamma_3),\quad \alpha\in(0,1)\\
&& \exists c_0\in\mathbb{R}\quad/ \quad \forall y \in \Omega \quad
:\qquad div(a(x)(x-y)) \leq c_0 \quad \mbox{ in }
\mathcal{D}^{\prime} (\Omega)\\
&&\Gamma_3 \text{ is } C_{loc}^{1,\alpha}
\end{eqnarray}

\vs0.2cm\n Here is the main result of this paper:

\begin{theorem}\label{t3.1} Let $w_0\in \pi_{x_1}(\Omega\cap\{x_2=h\})$ such that
 $(w_0,\Phi_h(w_0))\in D_h$, $T_h(\alpha_+(w_0),w_0))\in \Gamma_3$ and
\begin{equation}\label{e3.5}
\left[\frac{|Y_h|\beta(x,\varphi)}{H.\nu}\right](X(\alpha_+(w_0),w_0)<Y_h(X(w_0,\Phi(w_0)))
\end{equation}
Then $\Phi_h$ is continuous at $w_0$.
\end{theorem}

\vs0,3cm\n\emph{Proof.} Let $w_0\in \pi_{x_1}(\Omega\cap\{x_2=h\})$ as in the theorem.
Since $T_h(\alpha_+(w),w))$ is continuous at $w_0$ and $\Gamma_3$ is relatively open
in $\partial\Omega$, there exists  $w_1<w_0$ and $w_2>w_0$ such that
\[T_h(\alpha_+(w),w))\in \Gamma_3\quad\text{for all } w\in(w_1,w_2)\]

\n From Lemma 2.2, we know that $(\widetilde{u},\widetilde{\chi})$ is a solution
on the domain
\[Z=\{(t,w):~w_1<w<w_2 ~\text{ and }~ h< t < \alpha_{+}(w)\}\]
of a similar problem to $(P_0)$. Therefore it is enough to check that
the assumptions of Theorem 4.1 of \cite{[S]} are satisfied.

\n First, we deduce from Proposition 1.2 that the function $\texttt{h}$ satisfies
$$\left\{%
\begin{array}{ll}
   0< \underline{h}\leq \texttt{h}(t,\omega) \leq C \bar{h} & \hbox{for a.e }(t,\omega)\in D_h \\
    0\leq \texttt{h}_t (t,\omega)  \leq C \bar{h} & \hbox{for a.e }(t,\omega)\in D_h .\\
\end{array}%
\right.$$

\n Next, since $H\in C^{1,1}_{loc}( \Omega)$, it is easy to see that
$\texttt{a}\in C^{0,1}( D_h )$. Then by arguing as in \cite{[ChaL3]}, we can show that
we have for some positive constant $c_0, C_0$
\begin{eqnarray*}
&&|\texttt{a}(t,\omega)|  \leq C_0\\
&&\texttt{a}(t,\omega)\xi.\xi\geq  c_0 |Y_h|  \xi|^2 \geq  c_0|\xi|^2
\qquad \forall(t,w)\in D_h~~\forall \xi \in \mathbb{R}^2
\end{eqnarray*}
\qed

\vs0,3cm\n Moreover, since we have on $\widetilde{\Gamma_3}$
\begin{eqnarray*}
\lambda(.,\widetilde{\varphi})-\texttt{h}\nu_2&=&{{|Y_h|}\over
{\sqrt{1+\alpha_{+}^{\prime 2}(w)}}}.{{\beta(., \varphi)(T_h(\alpha_{+}(w),w)) }\over
{H.\nu(T_h(\alpha_{+}(w),w))}}-|Y_h|(\alpha_{+}(w),w)\nu_2\\
&=& |Y_h|\Big[{{\beta(., \varphi)}\over
{H.\nu}}-1\Big](T_h(\alpha_{+}(w),w))\nu_2
\end{eqnarray*}

\n this function is continuous on $\widetilde{\Gamma_3}$.

\n Finally, arguing as in the proof of Theorem 2.1 and using (3.5), we can show that
\begin{eqnarray*}
{{\lambda((\alpha_+(w_0),w_0),\widetilde{\varphi}(\alpha_+(w_0),w_0))}
\over{\texttt{h}(\phi_h(w_0),w_0)\nu_2(\alpha_+(w_0),w_0)}}
&=& {{|Y_h|\beta(., \varphi)(T_h(\alpha_{+}(w_0),w_0)))(\alpha_{+}(w_0),w_0)}\over
{|Y_h|(\phi_h(w_0),w_0)H.\nu(T_h(\alpha_{+}(w_0),w_0))}}<1\\
\end{eqnarray*}

\n We conclude that the function $\phi_h$ is continuous at $w_0$.

\qed

\end{document}